\documentclass[12pt]{article}

\usepackage[usenames]{color}
\usepackage{colortbl}
\usepackage{tikz-cd}
\usepackage{xcolor}
\usepackage[bbgreekl]{mathbbol}

\usepackage{amscd,amsmath, amssymb, fancyhdr, mathbbol,  url, euscript, amssymb}
\usepackage[backref=page]{hyperref}
\usepackage{mathtools}

\hypersetup{
	colorlinks   = true,
	citecolor    = magenta,
	linkcolor    =blue,
	urlcolor     =magenta	
}

\numberwithin{equation}{section}
\newcommand{\version}{}

\def\eqref#1{(\ref{#1})}

\def\C{{\Bbb C}}
\newcommand{\R}{{\Bbb R}}
\newcommand{\Q}{{\Bbb Q}}
\renewcommand{\H}{{\Bbb H}}

\def\1{\sqrt{-1}\:}

\newcommand{\cntrct}                
{\hspace{2pt}\raisebox{1pt}{\text{$\lrcorner$}}\hspace{2pt}}

\renewcommand{\tilde}{\widetilde}
\renewcommand{\bar}{\overline}
\renewcommand{\phi}{\varphi}
\renewcommand{\epsilon}{\varepsilon}
\renewcommand{\geq}{\geqslant}
\renewcommand{\leq}{\leqslant}

\newcommand{\be}{\begin{equation}}
\newcommand{\bea}{\begin{eqnarray}}
\newcommand{\eea}{\end{eqnarray}} 
\newcommand{\ee}{\end{equation}}

\newcommand{\End}{\operatorname{End}}

\newcommand{\Id}{\operatorname{Id}}

\newcommand{\Vol}{\operatorname{Vol}}

\newcommand{\Hom}{\operatorname{Hom}}

\newcommand{\codim}{\operatorname{codim}}

\newcommand{\Tr}{\operatorname{Tr}}

\newcommand{\vp}{\varphi}

\newcommand{\ve}{\varepsilon}
\newcommand{\p}{\partial}

\newcommand{\na}{\nabla}

\newcommand{\YM}{\text{YM}}
\newcommand{\HYM}{\text{HYM}}

\renewcommand{\Re}{\operatorname{Re}}

\renewcommand{\o}{\omega}

\newcounter{Mycounter}[section]
\newcounter{lemma}[section]
\newcounter{claim}[section]
\newcounter{sublemma}[section]
\newcounter{corollary}[section]
\newcounter{theorem}[section]
\newcounter{conjecture}[section]
\newcounter{proposition}[section]
\newcounter{definition}[section]
\newcounter{example}[section]
\newcounter{remark}[section]
\newcounter{problem}[section]
\newcounter{question}[section]
\makeatletter

\@addtoreset{equation}{section} \@addtoreset{footnote}{section}
\makeatother

\begin{document}
\title{Limiting properties of Hermitian-Einstein metrics under variation of the K\"ahler class\\ {\small Nikita Klemyatin}}
\date{}
\maketitle

\begin{abstract}
We study the behaviour of Hermitian-Einstein metrics on a holomorphic vector bundle over a compact K\"ahler manifold, as the underlying K\"ahler metric varies and its K\"ahler class changes. In particular, if the bundle stops being stable at some point, the Hermitian-Einstein metrics blow up, and their rescaling can be used to construct a destabilizing subsheaf, following the argument of Uhlenbeck and Yau. Furthermore, we study the behaviour of Chern connections of these metrics, and show, that after passing to a subsequence, and applying suitable gauge transforms these connections converge, and identify the blow-up locus of their Uhlenbeck limit.
	{\footnotesize }
\end{abstract}


\section{Introduction}

The relationship between the existence of Hermitian-Einstein metrics and the stability of holomorphic vector bundles has a long and profound history. Starting with the work of Narasimhan and Seshadri, this culminated in the works of Donaldson \cite{Donaldson85, Donaldson87} and Uhlenbeck-Yau \cite{UhlenbeckYau86}, where it was proved that the slope (poly)stability of a holomorphic vector bundle is a necessary and sufficient condition for the existence of the Hermitian-Einstein metric. This result was extended by many other authors, most notably by Simpson in \cite{Simpson} to Higgs bundles, and Bando and Siu \cite{BandoSiu} for stable sheaves on compact manifolds. Furthermore, Li and Yau \cite{LiYau} extended the original Donaldson-Uhlenbeck-Yau theorem to non-K\"ahler manifolds.

For non-stable bundles the research was started by Daskalopoulos and Wentworth \cite{DasWentworth1, DasWentworth2}, and then later expanded by Jacob \cite{Jacob1, Jacob2}, Sibley \cite{Sibley}, Collins-Jacob \cite{CollinsJacob} and Sibley-Wentworth \cite{SibleyWentworth}. Their combined result shows that outside an analytic subvariety 

\[Z_{alg} = \{x \in X ~|~\text{Gr}^{hns}(E)_x ~\text{is not free}\},\] 

there is a Yang-Mills connection on $E|_{X \setminus Z_{alg}}$.

However, the notion of the slope stability depends on the metric on the base manifold $X$. Moreover, it is known that a bundle $E$ may be stable with respect to one metric $\o_1$, but fail to be stable with respect to another one $\o_0$. It is natural to ask what happens with the corresponding Hermitian-Einstein metrics, if one deforms the underlying metric on $X$ in an unstable direction. This question has natural parallels with convergence of Einstein metrics on Riemannian manifolds (see \cite{Andersen1,Andersen2,CheegerNaber}), and cscK metrics on compact K\"ahler manifolds  when the K\"ahler class approaches the canonical class (see \cite{GuoJianShiSong,ChenWeber}). 

Another way to look at the question is through the variation of symplectic (or GIT) quotients. Let $(X, \o)$ be a compact K\"ahler manifold, and let $(E, H)$ be a hermitian vector bundle over $X$. Let $A$ be an integrable Hermitian connection on $E$. The set of integrable connections $\mathcal{A}^{1,1}_H$ can be written as $A + a$, where $a \in \Gamma(X, \mathfrak{u}(E)\otimes \Lambda^1_X)$, and $a$ satisfies the Maurer-Cartan equation $\bar{\partial}_A a^{0,1} + a^{0,1} \wedge a^{0,1} = 0$. One can define a symplectic form on $\mathcal{A}^{1,1}_H$ as:
\[\Omega(a,b) = -\int_X \Tr(a \wedge b) \wedge {\o^{n-1} \over (n-1)!}.\]

The complex gauge group $\mathcal{G}_\C$ acts on $\mathcal{A}^{1,1}_H$ in the standard way, and the real group $\mathcal{G}_\R(E)$, which consists of the isometries of $(E,H)$, preserves the symplectic form. In this setting, one can define the corresponding moment map, and show, that the zeroes of this moment map in a given orbit of $\bar{\p}_{A}$ are the Hermitian Yang-Mills connections, i.e. the connections, that satisfy the equation
\[\sqrt{-1}\Lambda_\o F_A - \mu_\o(E)\Id=0.\] 
If such $A$ exists, then it is standard that one can find a Hermitian-Einstein metric on $E$.

Now if we vary the symplectic form, then the corresponding moment map changes, and at some point, it may happen, that the moment map has no zeroes. A similar variation question in the case of actions of reductive groups on projective manifolds was addressed in \cite{DolgachevHi}, and in our framework these questions were considered in \cite{GrebToma, Zqin}. In these works, the authors studied wall-crossing and stability chambers, mostly on projective manifolds. Also, the related questions for the Gieseker stable bundles were considered in \cite{GrebRossToma,MatsukiWentworth}. However, the differential-geometric side of the picture was much less studied.  

However, it is easy to see, that we must impose additional assumptions on the Hermitian-Einstein metrics, as one can see, that for different K\"ahler metrics in a given class $[\o]$ the corresponding Hermitian-Einstein metrics may be different as well. Also, one must suitably normalize Hermitian-Einstein metrics, as multiplication by a positive number does not change the curvature of the metric, while it clearly may affect convergence. 

Throughout this paper we work under the following set of assumptions:

\hfill

\textbf{Assumptions:}

\begin{enumerate}\label{Assumptions}
	\item  $(X, \o_0)$ is a compact K\"ahler manifold of complex dimension $n$;
	
	\item $\o_s,~s \in (0;1]$ is a smooth family of uniformly equivalent K\"ahler metrics. This means that each $\o_s$ is smooth, they depend smoothly on $s$, and there is a uniform constant $C>1$, such that for all $s \in [0;1]$ we have
	
	\[C^{-1}\o_0 \leq \o_s \leq C\o_0;\]
	
	We also require that $\o_s$ tends to $\o_0$ uniformly as $s \to 0$. Later, in Section \ref{UhlenbeckCompSection}, we also need the additional assumption that curvatures of $\o_s$ are uniformly bounded.
	
	\item $E$ is a holomorphic vector bundle of rank $r$, which is stable with respect to all $\o_s, ~s \in (0;1]$;
	
	\item Let $H_s$ be a corresponding family of Hermitian-Einstein metrics (HE), depending smoothly on $s \in (0;1]$ (this should be possible to arrange, via the implicit function theorem, \cite[iv]{Siu}). More precisely:
	
	\[\sqrt{-1}\Lambda_{\o_s}F_{H_s} = \mu_{\o_s}(E) \Id,\]
	
	where
	
	\[\mu_{\o_s}(E) = {2\pi \over r \Vol(X, \o_s)} \int_X c_1(E) \wedge {\o^{n-1}_s \over (n-1)!}\]
	
	\item We also introduce the relative endomorphism $h_s = H^{-1}_1H_s$, and we can choose $H_s$ such that $\int_X \log \det h_s \o^n_s = 0$. We will show that under this assumption, the determinant of $h_s$ is uniformly bounded above and below.
	
\end{enumerate}

A few comments about our assumptions. As we discussed earlier, the definition of $\mu_\o(E)$, and hence, the stability notion, depend on the K\"ahler class rather than the metric itself. However, if one changes the metric within the K\"ahler class, then the corresponding Hermitian-Einstein metric also changes. If we have a smooth family $[\o_s]$ of K\"ahler classes, then the corresponding family of K\"ahler metrics, and the corresponding family of Hermitian-Einstein metrics $H_s$ need not be smooth. Hence, it is necessary to vary the K\"ahler metrics themselves. Moreover, if we want to study compactness properties of the Hermitian-Einstein metrics $H_s$, it is necessary to impose a certain normalization on these metrics, and we later show that our normalization for the relative endomorphism $h_s$ is natural to assume, and can always be achieved. It is not hard to see that our assumptions on K\"ahler metrics are quite natural too, and do not impose any additional restrictions on their K\"ahler classes. Indeed, if $[\o_0]$ lies on the boundary of a stability chamber, and $[\o_1]$ is the class, lying within this chamber and sufficiently close to $[\o_0]$, then the family  $\o_s = (1-s)\o_0 + s\o_1$ satisfies all our regularity assumptions. One can construct other examples of such families of K\"ahler metrics.

Our main result here is the following.

\hfill

\begin{theorem}\label{MainTheorem}
	Let $X$ be a compact complex manifold, and suppose that the Assumptions above hold. Let $H_s$ be the corresponding family of Hermitian-Einstein metrics, such that the endomorphism $h_s = H^{-1}_1H_s$ is normalized by $\int_X \log\det h_s \o^n_s =0$. Then
	\begin{enumerate}
		\item If $\Tr h_s$ is uniformly bounded, then $E$ is $\o_0$-stable.
		
		\item If $\Tr h_s$ is unbounded, then the bundle $E$ is strictly $\o_0$-semistable. The projector $\pi_{\mathcal{F}}$ onto a $\o_0$-destabilizing sheaf $\mathcal{F}$ can be constructed as a weak limit of endomorphisms $\Id - \tilde{h}_s^\sigma = \Id - \Big({h_s \over \sup_X \Tr h_s}\Big)^\sigma$.
				
		\item If we have a uniform bound on $\text{Rm}(\o_s)$, then there is a subsequence $s_j$, and a subset $Z_{an}$, such that the Chern connections weakly converge (after a suitable gauge transform) in $L^p_{1,loc}(X \setminus Z_{an}), ~p>2n$, and blow up along $Z_{an}$. If, in addition, we have $C^k$ bounds on $\text{Rm}(\o_s)$, then the convergence of connections can be taken in $C^k_{loc}(X \setminus Z_{an})$. Also, the limiting bundle $E_\infty$ can be extended over all of $X$ as a reflexive coherent sheaf, and we have equality
		\[E_\infty \cong (\text{Gr}_{\o_0}^S(E))^{**}.\]
		
		\item In addition, under the same assumptions as previously, the set  $Z_{an}$ is an analytic subvariety, and the sets $Z_{an}$ and $Z_{alg}$ coincide. Any irreducible component $Z \subset Z_{an}$ of codimension 2 lies in $Z_{alg}$, and the multiplicities of $Z$ in $Z_{an}$ and $Z_{alg}$ are the same.
	\end{enumerate} 
\end{theorem}

\hfill

The theorem says, that, as we approach the wall of the stability chamber, the appropriate choice of Hermitian-Einstein metrics for $[\o_s]$ detects instability of our bundle $E$ with respect to $[\o_0]$. Moreover, the destabilizing sheaves and the Seshadri filtration emerge directly from these metrics. 

\ref{MainTheorem} should be compared with the results from \cite{ClarkeTipler}. They establish results similar to part (3) of \ref{MainTheorem}, under the additional assumption that $E$ is sufficiently smooth, i.e. $\text{Gr}^S_{\o_0}(E)$ is locally free, and has abelian automorphism group. Our result does not require such restrictions. 

It is also worth mentioning, that the results (1)-(2) of \ref{MainTheorem} hold if our metrics $\o_s$ are balanced or Gauduchon. The Uhlenbeck compactness from (3) also can be proved, under additional assumptions on the curvatures and torsion of $\o_s$. The part (4) of the theorem cannot be proven for non-K\"ahler manifolds via our methods, as it requires certain homological properties of compact K\"ahler manifolds. The paper \cite{Delloque} also gives results similar to part (3) of \ref{MainTheorem} in the setting of balanced metrics, and also removes the abelian automorphism assumption. However, the assumption on sufficient smoothness remains in that paper.

The organization of the paper is as follows: in Section 2 we review all the necessary background material, set up the notation, and establish the normalization for $h_s$. We also establish in \ref{MainReason} the first part of \ref{MainTheorem}. In Section 3 we adapt the method of Uhlenbeck and Yau from \cite{UhlenbeckYau86}, and prove the second part of \ref{MainTheorem}. Section 4.1 is devoted to preliminaries for Uhlenbeck compactness for Hermitian Yang-Mills connections in the case of varying metrics $\o_s$. In Section 4.2 we establish Uhlenbeck compactness in our setting, and show part (3) of \ref{MainTheorem}. In Section 4.3 we show that $Z_{an}$ is rectifiable and analytic. Finally, in Section 4.4 we establish the last part of \ref{MainTheorem}.

\textbf{Acknoledgements:} The author is grateful to Prof. Yingying Zhang and to Prof. Shing-Tung Yau for suggesting the problem, and their encouragement. 

\textbf{AI disclosure:} The author used AI to check the language, exposition, literature review, and to serve as a proofreader. During this review, AI found a gap in the initial proof of \ref{MainReason},  as well as several other minor issues throughout the text. The correction of these issues was done by the author. All other proofs were also done by the author, and the author bears exclusive responsibility for the correctness of the results.   

\section{Preliminaries.}\label{PrelimSection}

\subsection{Harder-Narasimhan and Seshadri filtrations.}

Let $(X,\o)$ be a compact K\"ahler manifold with volume $V = \int_X {\o^n \over n!}$, and $E \to X$ be a holomorphic vector bundle over $X$ of rank $r$. We define {\it the slope} of $E$  by the formula

\begin{equation}\label{Slope}
	\mu_\o(E) = {2 \pi \over r V}\int_X c_1(E) \wedge {\o^{n-1} \over (n-1)!}.
\end{equation}

This definition is also valid for coherent sheaves over $X$.

The bundle is called (semi)stable if for any proper and nonzero torsion-free coherent subsheaf $\mathcal{F}$ of $E$ with $0 < \text{rank}(\mathcal{F}) < r$ there is an inequality

\[\mu(\mathcal{F}) (\leq) <\mu(E).\]

Clearly the definition of stability depends on the K\"ahler class, rather than the K\"ahler form itself. Moreover, if $E$ is stable with respect to $[\o]$, then it is also stable with respect to any $[\tilde{\o}]$ sufficiently close to $[\o]$. Hence, the set 
\[\text{Stab}(E) = \{ [\o] \in H^{1,1}(X, \R) ~|~ [\o] ~\text{is K\"ahler and}~ E ~\text{is}~ [\o]\text{-stable} \}\]

is open and stable under dilations $[\o] \to \lambda [\o], ~\lambda>0$.

Next, we recall the notion of Harder-Narasimhan (HN) filtration.

\hfill

\begin{definition}\label{HNFiltration}(see \cite[Theorem 5.7.15]{Kobayashi}, \cite[Definition Proposition 2.2]{DasWentworth1})	
	Given any torsion free coherent sheaf $\mathcal{E}$ over $(X, \o)$, there exists a filtration 
	\[0 = \mathcal{E}_0 \subset \mathcal{E}_1 \subset \dots \subset \mathcal{E}_m = \mathcal{E},\]
	such that each $\mathcal{Q}_{j+1} = \mathcal{E}_{j+1} / \mathcal{E}_j$ is semistable and torsion-free with $\mu_{\o}(\mathcal{Q}_j) > \mu_{\o}({\mathcal{Q}_{j+1}})$. Such a filtration is called the {\it Harder-Narasimhan} filtration of $\mathcal{E}$. The associated graded object $\text{Gr}^{HN}(\mathcal{E}) = \oplus_{j=1}^m \mathcal{Q}_j$ is uniquely defined by the isomorphism class of $\mathcal{E}$.
\end{definition}

\hfill

For semistable sheaves there is another filtration, which one may consider as a decomposition of such sheaves into stable parts.

\hfill

\begin{definition}\label{SFiltration}(see \cite[Proposition 2.5]{DasWentworth1})
	Given a semistable sheaf $\mathcal{Q}$ over $(X, \o)$, there is a filtration
	\[0 = \mathcal{F}_0 \subset \dots \subset \mathcal{F}_l = \mathcal{Q},\]
	such that for any $j \in \{1, \dots, l\}$ the sheaf $\mathcal{F}_{j } /\mathcal{F}_{j-1}$ is stable and torsion-free, and $\mu_\o (\mathcal{F}_{j } /\mathcal{F}_{j-1}) = \mu_\o(\mathcal{Q})$. Such a filtration is called {\it a Seshadri filtration}. The double dual of the associated graded object $\text{Gr}^S(\mathcal{Q}) = \oplus_j \mathcal{F}_{j } /\mathcal{F}_{j-1}$ is uniquely determined by the isomorphism class of $\mathcal{Q}$ (see \cite{ChenSun}).
\end{definition}

\hfill

For a semistable sheaf $\mathcal{Q}$ one can define the notion of algebraic singular set.

\hfill

\begin{definition}\label{ZalgDefinition}
	Let $\mathcal{Q}$ be a semistable torsion-free sheaf over a compact K\"ahler manifold $(X, \o)$. The algebraic singular set is defined as
	\[Z_{alg} = Z_{alg}(\mathcal{Q}) = \{ x \in X ~|~ \text{Gr}^S(\mathcal{Q})_x ~\text{is not free}\}.\]
	It is an analytic subvariety of complex codimension at least two. 
\end{definition}

\hfill

We conclude this section by an observation: as we pointed out before, the notions of stability, as well as Harder-Narasimhan and Seshadri filtrations depend on the background K\"ahler class $[\o]$. With a different choice of $[\o]$ the stability may change, as well as the filtrations above. There are examples (see \cite[p.99]{Friedman}), that already on algebraic surfaces there exist bundles, which are stable with respect to one polarization, and not stable with respect to another.

\subsection{Hermitian-Einstein metrics, stability and the Einstein constant.}
Let $(X,\o)$ be a compact K\"ahler manifold, and $E$ be a holomorphic bundle over $X$ with a Hermitian metric $H$. Recall that the metric is called Hermitian-Einstein, if
\[\sqrt{-1}\Lambda_\o F_H = \lambda \Id.\]

A direct computation shows that $\lambda = \mu_\o(E)$ in our case.

There is the following fundamental theorem, which relates the existence of Hermitian-Einstein metric to the questions of stability.

\hfill

\begin{theorem}\label{DUY}
	The bundle $E$ over $(X, \o)$ admits a Hermitian-Einstein metric if and only if it is polystable, i.e. it is a direct sum of stable vector bundles with the same slope.
\end{theorem}

\hfill

The proof can be found in \cite{Donaldson85, Donaldson87} for the case of projective manifolds, and in \cite{UhlenbeckYau86, Simpson} for general K\"ahler manifolds.

We now briefly discuss the formalism, that leads to one direction of \ref{DUY}. Assume that we have a short exact sequence of bundles over $X$:
\[0 \to S \to E \to Q \to 0.\]

The metric $H$ induces metrics on $S$ and $Q$ in the usual way.

The curvatures are related by the following identity:
\begin{equation}
	\1 F_E = 
	\begin{pmatrix}
		\1 F_S  + \1 \gamma \wedge \beta & \1  \p_\nabla \gamma \\
		\1  \bar{\p}\beta & \1 F_Q + \1 \beta \wedge \gamma
	\end{pmatrix}
\end{equation}

Here $\beta = -\gamma^\dagger =  \na \pi_S$ is the second fundamental form of $S$. The $(0,1)$-form $\gamma$ lies in $\Gamma(X, \Hom(Q, S) \otimes \Lambda^{0,1}_X)$, and satisfies the equation $\bar{\p}\gamma = 0$.

Since the formula above is local, we obtain the following useful statement. 

\hfill

\begin{proposition}(Chern-Weil formula, see \cite[Lemma 3.2]{Simpson})\label{ChernWeil}
	
	Let $\mathcal{F}$ be a saturated torsion-free subsheaf  of a holomorphic bundle $E$, and let $\pi=\pi_{\mathcal{F}} \in L^2_1(X,\o)$ be an orthogonal projector onto $\mathcal{F}$. We have the following formula for the slope $\mu_\o (\mathcal{F})$:
	
	\[\mu_\o(\mathcal{F}) = {1 \over \text{rank}(\mathcal{F})V_\o}\int_X \big(\1\Tr((\Lambda_\o F_E)\pi) - |\na \pi|^2\big) {\o^n \over n!}. \]
	
	Moreover, the integral on the right hand side may be computed over $X \setminus Z$, where $Z$ is a codimension at least two analytic set, where $\mathcal{F} \hookrightarrow E$ is not a subbundle. 
\end{proposition}

\hfill 

The Chern-Weil formula also implies the "easy" direction in \ref{DUY}.

We now return to our setting. Assume $(X, \o_s, E, H_s)$ as in the statement of \ref{MainTheorem}. We will denote $H_1$ by $H$, and use it as a given reference metric. We have the following important proposition.

\hfill

\begin{proposition}\label{NormalizationOfEnd}
	Consider the relative endomorphism $h_s = H^{-1}H_s$, normalized by
	\[\int_X \log \det h_s \o^n_s = 0.\]
	Then $\det h_s$ is uniformly bounded, i.e. there exists $C>1$, such that
	\[C^{-1} \leq \det h_s \leq C.\]
\end{proposition}

\hfill

\begin{proof}
	Note that
	\[\Delta_s \log \det h_s = \sqrt{-1}\Tr\Lambda_{\o_s}F_H - r\mu_{\o_s}(E) .\]
	By our assumption, the RHS is uniformly bounded. By the result of \cite{GPSSGreen}, the Green functions of the metrics $\o_s$ are uniformly bounded from below, with the bound independent of $s$. Let $u_s = \pm \log \det h_s$. The Green formula
	\[u_s = {1 \over \Vol(X, \o_s)}\int_X u_s \o^n_s + \int_X G_{\o_s}(x,y)(-\Delta_s u_s)(y)\o^n_s(y)\] 
	implies that $u_s \leq C$, where $C >0$ depends on the lower bound on the Green function, $\mu_{\o_s}(E)$ and $F_H$.
\end{proof}

\hfill

Note that in our estimate we only need the uniform equivalence of $\o_s$. Of course, one can do the usual Moser iteration argument, but the estimates from \cite{GPSSGreen} allow us to give a particularly streamlined proof in the K\"ahler case.

The next proposition is important, and deals with the first part of \ref{MainTheorem}. 

\hfill

\begin{proposition}\label{MainReason}
	
	If the endomorphisms $h_s = H^{-1}H_s$ have uniformly bounded eigenvalues for $s \in (0;1]$, then there is a subsequence $\{s_j\}$, such that $H_j = H_{s_j} \to H_0$ as $s_j \to 0$, and $H_0$ is a Hermitian-Einstein metric. Moreover, $E$ is $\o_0$-stable.
\end{proposition}

\hfill

\begin{proof}
	
	\textbf{Step 1:}
	
	First of all, we need to show the following statement 
	
	{\it Suppose $E$ is a trivial holomorphic bundle over $\C^n$. Let $H$ be a hermitian metric on $E$, which solves $\Lambda_{\o_{euc}}F_{H} = 0$ on $\C^n$, and $h = H_{flat}^{-1}H$ is the relative endomorphism with respect to the standard metric on $E$. If $h$ has uniformly bounded eigenvalues, then $H$ is constant.}
	
	Let $S = \na h h^{-1}$, where $\na$ is the Chern connection of $H$. The standard formula, together with upper and lower bounds on $h$ gives:
	\[\Delta \Tr h \geq c_0 |\na h h^{-1}|^2.\]
	
	Here the norm is taken with respect to $H$. 
	
	Then
	\[\Delta |S|^2 = |\na S|^2 + |\bar{\na}S|^2 + 2\Re \na^{\bar k}\na_{\bar k} S \cdot S \geq 0,\]
	because flatness of $\C^n$ together with $\Lambda_{\o_{euc}}F_H = 0$ imply $\na^{\bar k}\na_{\bar k} S = \na_{\bar k} \na^{\bar k}S = 0$.
	
	Thus, $|S|^2$ is subharmonic.
	
	Using the standard cutoff function $\chi_R, ~\chi_R \equiv 1 ~\text{on}~ B_R$, supported on $B_{2R} \subset \C^n$, and $|\Delta \chi_R| \leq CR^{-2}$, we obtain
	\[\int_{B_R}|S|^2 \leq {1 \over c_0}\int_{B_{2R}} \chi_R \Delta \Tr h\leq CR^{2n-2}   . \]
	
	Hence, by the mean value inequality for subharmonic functions, we obtain
	\[|S|^2(0) \leq {C \over R^{2n}}\int_{B_R(0)}|S|^2 \leq CR^{-2}.\]	
	
	Hence $|S|^2(0) = 0$. Since the same argument applies to any point, this means $\na h = 0$, hence $H$ is constant.
	
	\textbf{Step 2:}
	The metrics $H_s$ satisfy 
	\[\1 \Lambda_{\o_s}F_{H_s} = \mu_s(E)\Id_E,\]
	and with our assumptions, the eigenvalues of $H_s$ with respect to $H=H_1$ are uniformly bounded above and below.
	
	Then we may adapt the argument from \cite[Lemma 19]{Donaldson85} or \cite{Siu}, to obtain the uniform bound on $\na^H h_s$. Suppose the gradients are unbounded. Since the metrics $\o_s$ converge to $\o_0$ uniformly, we can find a sequence $h_j$, and corresponding sequence of metrics $\o_j$, such that 
	\[ \sup_X |\na^H h_j| = |\na^H h_j|(x_j) = M_j \to + \infty, ~ x_j \to x_0. \]
	 We rescale the metrics in the usual way
	\begin{align*}
		\hat{h}_j(x) = h_j(x_j + {x \over M_j}), ~ \hat{\o}_j(x) = \o_j(x_j + {x \over M_j}).
	\end{align*}
	
	Note that $\hat{\o}_j$ converge uniformly to $\o_{euc}$.
	
	Now, the uniform equivalence and uniform convergence of the metrics $\o_s$ allow us to obtain uniform $L^p_2$-estimates for $\hat{h}_j$ (see \cite[Theorem 9.11]{GilbargTrudinger}). Then $\hat{h}_j$ (up to taking a subsequence) converge in $C^1$ to some metric $\hat{H}_0$, which satisfies 
	\[\Lambda_{\o_{euc}}F_{\hat{H}_0} = 0. \]
	
	Step 1 guarantees that $\hat{H}_0$ is constant. Hence $\na \hat{H}_0 = 0$, which contradicts the fact that $|\na^H \hat{h}_j|(0) = 1$. Hence, the original family has bounded gradients, and there is a subsequence $H_j$, which converges to $H_0$ in $C^1$.
	
	Let $H_0 = H h_0$. Then $H_0$ satisfies $\sqrt{-1}\Lambda_{\o_0}F_{H_0} = \mu_{\o_0}(E)\Id$ weakly. Now the standard elliptic bootstrap implies that $H_0$ is smooth and Hermitian-Einstein. Clearly $E$ is $\o_0$-stable, because it is stable with respect to other metrics $\o_s$. 
\end{proof}

\hfill

\begin{remark}
	We could avoid Step 1 in the proposition above by assuming stronger control of the metrics $\o_s$.
	
	Let $S = \na h_s h^{-1}_s$, where $\na$ is the Chern connection with respect to $H_s$. A straightforward calculation for $|S|^2 = |S|_{H_s}^2$ yields:
	\[\Delta_s |S|^2 \geq -C_1(|S|^2 + 1),\]
	where $C_1>0$ depends on $H_1$, and the bounds on $\o_s$ and $\text{Ric}(\o_s)$. On the other hand, if $\Tr h_s$ is uniformly bounded, then there are uniform constants $C_2, C_3$, such that
	\[\Delta_s \Tr h_s \geq -C_2 + C_3|S|^2. \]
	
	Now we may apply the maximum principle to $|S|^2 + A\Tr h_s$, and show that $|S|^2$ is uniformly bounded. However, this argument requires control of the Ricci curvature of $\o_s$.
\end{remark}

\subsection{Yang-Mills, Hermitian-Yang-Mills functionals, and their generalizations.}

Let $A$ be an integrable connection on a hermitian bundle $(E,H)$ over a compact K\"ahler manifold. We define the Yang-Mills (YM) functional as
\[\YM_\o(A):= \int_X |F_A|^2_H {\o^n \over n!}.\]

The Hermitian Yang-Mills (HYM) functional is defined as follows:
\[\HYM_\o(A) := \int_X |\Lambda_\o F_A|^2_H {\o^n \over n!}.\]

Hence $\YM_\o$ and $\HYM_\o$ have the same critical points. 
They are related by the following identity (see \cite{Kobayashi}):
\begin{equation}\label{YMHYM}
	\YM_\o(A) = \HYM_\o(A) + {4\pi^2 \over (n-2)!}\int_X (2c_2(E) - c_1^2(E)) \wedge \o^{n-2}.
\end{equation}

In fact, we can reformulate this as follows:
\begin{align}\label{YMHYM2}
	&\YM_\o(A) = ||\sqrt{-1}\Lambda_\o F_A - \mu_\o(E)\Id||^2_{L^2(X,\o)} + r \mu_{\o}(E)^2 V_{\o}
	\\
	&+ {4\pi^2 \over (n-2)!}\int_X (2c_2(E) - c_1^2(E)) \wedge \o^{n-2}.
\end{align}

Moreover, critical points of $\YM_\o$ on the set $\mathcal{A}^{1,1}_H$ of integrable hermitian connections are generalized Hermitian Yang-Mills connections (see \cite{Kobayashi}).  

In the next sections we will need the following useful lemma.

\hfill

\begin{proposition}\label{InfAttained}
	Let $\o_s$ be a family of K\"ahler metrics, that satisfy our assumptions, and let $H_s$ be the family of the Hermitian-Einstein metrics on $E$ with respect to $\o_s$. Then
	\begin{gather*}
		\lim_{s \to 0}\YM_{\o_s}(A_{H_s}) = \lim_{s \to 0}\YM_{\o_0}(A_{H_s})
		\\
		= r\mu^2_{\o_0}(E)V_{\o_0} + {4\pi^2 \over (n-2)!}\int_X (2c_2(E) - c_1^2(E)) \wedge {\o_0}^{n-2},
	\end{gather*}
	and
	\[\lim_{s \to 0}\HYM_{\o_s}(A_{H_s}) = \lim_{s \to 0}\HYM_{\o_0}(A_{H_s}) = r\mu^2_{\o_0}(E)V_{\o_0}.\]
	
	In particular, $\sqrt{-1}\Lambda_{\o_0} F_{A_{H_s}} \xrightarrow[s \to 0]{L^2} \mu_{\o_0}(E)\Id$, and 
	\[\lim_{s \to 0}\YM_{\o_s}(A_{H_s}) = \inf_{A \in \mathcal{A}^{1,1}_H} \YM_{\o_0}(A).\] 
\end{proposition}

\hfill

\begin{proof}
	By \ref{YMHYM} we have
	\[\YM_{\o_s}(A_{H_s}) = r\mu^2_{\o_s}(E)V_{\o_s} + {4\pi^2 \over (n-2)!}\int_X (2c_2(E) - c_1^2(E)) \wedge \o_s^{n-2},\]
	and
	\[\YM_{\o_0}(A) = \HYM_{\o_0}(A) + {4\pi^2 \over (n-2)!}\int_X (2c_2(E) - c_1^2(E)) \wedge \o_0^{n-2}.\]
	
	The integral $\int_X (2c_2(E) - c_1^2(E)) \wedge \o_s^{n-2}$ and $\mu_{\o_s}(E)$ are both continuous in $s$. Next, $\YM_{\o_s}(A_{H_s})$ is uniformly equivalent to $\YM_{\o_0}(A_{H_s})$. More precisely, there are $\ve_s>0, ~\ve_s \to 0$, such that
	\[(1 - \ve_s)\YM_{\o_s}(A_{H_s}) \leq \YM_{\o_0}(A_{H_s}) \leq (1+ \ve_s)\YM_{\o_s}(A_{H_s}).\] 
	
	Now, the statement follows from the fact that  
	\[\lim_{s \to 0}\YM_{\o_s}(A_{H_s}) = r\mu^2_{\o_0}(E)V_{\o_0} + {4\pi^2 \over (n-2)!}\int_X (2c_2(E) - c_1^2(E)) \wedge {\o_0}^{n-2},\]
	and \ref{YMHYM2}. 
	
	Alternatively, 
	\[||\sqrt{-1}\Lambda_{\o_0}F_{H_s} - \mu_{\o_0}(E)\Id ||_{L^2} \leq C ||\o_s - \o_0||_{C^0} ||F_{H_s}||_{L^2} + C|\mu_{\o_0}(E) - \mu_{\o_s}(E)|,\]
	which tends to zero as $s \to 0$. Hence 
	\[\sqrt{-1}\Lambda_{\o_0} F_{A_{H_s}} \xrightarrow[s \to 0]{L^2} \mu_{\o_0}(E)\Id.\] 
	The last statement immediately follows from this observation. 
\end{proof}

\hfill

Finally, we recall the definition of generalized Hermitian-Yang-Mills functionals, introduced in \cite{AtiyahBott, DasWentworth1}.

\hfill

\begin{definition}\label{GeneralizedHYM}
	Let $\alpha \in [1; +\infty)$, and $\Phi_\alpha(\mu) = \sum_{\gamma = 1}^r|\mu_\gamma|^\alpha$. We define $\HYM_{\alpha,\o}(A)$ as
	\[\HYM_{\alpha,\o}(A) := \int_X \Phi_\alpha(\sqrt{-1}\Lambda_\o F_A) {\o^n \over n!}.\]
	For any $N \in \R$ we similarly define $\HYM_{\alpha, N, \o}(A)$ as 
	\[\HYM_{\alpha, N, \o}(A) := \int_X \Phi_\alpha(\sqrt{-1}\Lambda_\o F_A + N\Id) {\o^n \over n!}.\]
\end{definition}

These functionals clearly behave like $L^\alpha$-norm of $\Lambda_\o F_A$. The main use of $\HYM_{\alpha, N, \o}$ functionals is the identification of the Harder-Narasimhan-Seshadri type of $A$. 

\section{Reduction to the method of Uhlenbeck and Yau}\label{ConstructionOfAProjector}

In this section we adapt the arguments of Uhlenbeck and Yau and construct a $\o_0$-destabilizing subsheaf of $E$ with the slope equal to the slope of $E$.  

We start with the standard estimate:

\hfill

\begin{lemma}\label{C2Estimate}
	Under the conditions above there is a constant $C>0$, independent of $s$, such that the following holds:
	\[\Delta_s \log \Tr h_s \geq - |\Lambda_s F_{H_s}| - |\Lambda_s F_1| \geq -C\]
\end{lemma}

\hfill

\begin{proof}
	The first inequality is a calculation in Siu's book \cite[1.9.2]{Siu}. The second one follows from the fact that the metrics $\o_s$ are comparable, and $H_s$ are Hermitian-Einstein. 
\end{proof}

\hfill

We immediately obtain (via Moser iteration), that there is $C>0$, independent of $s$, such that
\begin{equation}\label{MoserEstimate}
		\sup_X \Tr h_s \leq C \int_X \Tr h_s {\o^n_0 \over n!}.
\end{equation}

We rescale our relative endomorphism $h_s$:

\[\tilde{h}_s = \frac{h_s}{\sup_X \Tr h_s}.\]

\hfill

\begin{proposition}\label{SobolevEstimate}
	There is a constant $C>0$, such that for every $s \in (0;1]$ and any $\sigma \in (0;1]$ the endomorphisms $\tilde{h}^\sigma_s$ are uniformly bounded by $C$ in $L^2_1(X, \End(E))$, i.e. 
	\[||\tilde{h}^\sigma_s||_{L^2_1(X, \o_0)} \leq C.\]
\end{proposition} 

\hfill

\begin{proof}

We recall a general identity relating the curvatures $F = F_{H_s}$ and $F_1$:
\begin{equation}\label{CurvatureRelation}
	F_{j \bar k} - (F_1)_{j \bar k}  = - \p_{\bar k}(h^{-1}_s\na_j h_s) 
\end{equation}

where $\na$ is (always) the Chern connection for $H:=H_1$ (we omit the bundle indices).

Taking the $H$-inner product with $h^\sigma_s, ~\sigma\in(0;1]$ gives:
\begin{equation}
	-g^{j \bar k}_s \langle \p_{\bar k}(h^{-1}_s\na_j h_s), h^\sigma_s \rangle = \langle \sqrt{-1}\Lambda_{\o_s} F_s - \mu_{\o_s}\Id, h^\sigma_s\rangle - \langle \sqrt{-1}\Lambda_{\o_s} F_1 - \mu_{\o_s}\Id, h^\sigma_s\rangle
\end{equation}

The first summand on the right hand side is zero. The left hand side can be rearranged as follows:
\begin{align}
	&-g^{j \bar k}_s \langle \p_{\bar k}(h^{-1}_s\na_j h_s), h^\sigma_s \rangle = -g^{j \bar k}_s \p_{\bar k}\langle (h^{-1}_s\na_j h_s), h^\sigma_s \rangle + g^{j \bar k}_s \langle (h^{-1}_s\na_j h_s), \na_{ k} h^\sigma_s\rangle 
	\\
	& = -g^{j \bar k}_s \p_{\bar k}\langle (\na_j h_s), h^{\sigma-1}_s \rangle + g^{j \bar k}_s \langle h^{-1}_s\na_j h_s, \na_{ k} h^\sigma_s\rangle 
	\\
	& \geq - {1 \over \sigma}\Delta_s \Tr h^\sigma_s + |h^{-\sigma/2}_s \na h^\sigma_s|^2_{\o_s},
\end{align}

and the last inequality follows from calculations of Uhlenbeck and Yau.

We obtain the following inequality:
\begin{equation}\label{TheKeyInequality}
	- {1 \over \sigma}\Delta_s \Tr h^\sigma_s + |h^{-\sigma/2}_s \na h^\sigma_s|^2_{\o_s} \leq - \langle \sqrt{-1}\Lambda_{\o_s} F_{H_1} - \mu_{\o_s}\Id, h^\sigma_s\rangle,
\end{equation}

and, by the rescaling,
\begin{equation}\label{RescaledKeyInequality}
	- {1 \over \sigma}\Delta_s \Tr \tilde{h}^\sigma_s + |\tilde{h}_s^{-\sigma/2}\na \tilde{h}^\sigma_s|^2_{\o_s} \leq - \langle \sqrt{-1}\Lambda_{\o_s} F_{H_1} - \mu_{\o_s}\Id, \tilde{h}^\sigma_s\rangle.	
\end{equation}

Integrating against $\o^n_s$, and using the equivalence of K\"ahler metrics together with the fact that the eigenvalues of $\tilde{h}_s$ are at most 1, we get the following inequality:

\begin{equation}
	\int_X |\na \tilde{h}^\sigma_s|^2 {\o^n_0 \over n!} \leq C,
\end{equation}

where $C>0$ depends only on $(X, \o_0, \o_s, H_1)$ (and note that $\o_0$ and $\o_s$ are equivalent). Q.E.D.

\end{proof}

\hfill

The result of \ref{SobolevEstimate} enables the argument of Uhlenbeck and Yau. Let $s_j \xrightarrow[j \to + \infty]{} 0$ be any sequence, such that $\sup_X \Tr h_{s_j} \to + \infty$. Let $h_j = h_{s_j}$, and $\o_j = \o_{s_j}$. We have the following:

\hfill

\begin{proposition}\label{LimitExistence}\label{WeakHolo}
	
	There is a nonnegative self-adjoint element $\H \in L^2_1(X, \End(E))$, such that for any rational $\sigma \in (0;1]$ the sequence (maybe after taking a subsequence) $\tilde{h}^\sigma_j$ converges to $\H^\sigma$ in $L^2$ (and a.e.), and $\na \tilde{h}^\sigma_j$ converges to $\na \H^\sigma$ in $L^2$ weakly. Moreover, there exists a decreasing sequence $\sigma \in \Q \cap (0;1]$, such that $\pi = \lim_{\sigma \to 0}(\Id - \H^\sigma)$. The limit $\pi$ lies in $L^2_1(X, \o_0)$, and satisfies the following identities
	\begin{align*}
		&\pi = \pi^2 = \pi^* 
		\\
		&(\Id - \pi)\bar{\p}\pi = 0 ~\text{a.e.}
	\end{align*}
	 and defines a proper coherent subsheaf $\mathcal{F}$ of $E$.
\end{proposition}

\hfill

\begin{proof}	
	This is a standard fact, but we include the proof for the reader's convenience. The author learned it from Professor Duong H. Phong. 
	
	The convergence statements and equality $\pi = \pi^2 = \pi^*$ are trivial. The difficulty is the second one. Note that $(\Id - \pi)\bar{\p}\pi = -(\pi \na(\Id- \pi))^*$, so we prove that the right hand side is zero.
	
	We show, that 
	\[\pi\na (\Id  - \pi) = \lim (\Id - \tilde{h}_j^s)\na \tilde{h}_j^\sigma.\]
	
	Here, of course, the equality is vague. We shall take three limits, and the order of limits matters.
	
	 First of all, from the H\"older inequality we have the weak convergence of $(\Id - \tilde{h}_j^s)\na \tilde{h}_j^\sigma$ to $(\Id - \H^s)\na \H^\sigma$. Indeed, for any $\psi \in L^2(X, \Lambda^{1,0}_X \otimes \End(E))$ we have	
	\[\langle (\Id - \tilde{h}_j^s)\na \tilde{h}_j^\sigma ,\psi\rangle_{L^2} = \langle (\Id - \H^s)\na \tilde{h}_j^\sigma ,\psi\rangle_{L^2}  + \langle (\H^s - \tilde{h}_j^s)\na \tilde{h}_j^\sigma ,\psi\rangle_{L^2},\]
	hence we have the weak convergence, by the dominated convergence theorem.
	
	Next, since every eigenvalue $\mu_\alpha$ of $\tilde{h}_j$ lies between $0$ and $1$, it satisfies the inequality
	\[1 - \mu^s_\alpha \leq \frac{2s}{2s + \sigma}\mu^{-\sigma/2}_\alpha, ~~~\alpha = 1, \dots, r.\]
	
	This implies
	\begin{equation}
		\int_X |(\Id - \tilde{h}_j^s)\na \tilde{h}_j^\sigma|^2 \frac{\o^n_0}{n!} \leq \Bigg(\frac{2s}{2s + \sigma}\Bigg)^2\int_X |\tilde{h}^{-\sigma/2}_j\na \tilde{h}^{\sigma}_j|^2 \frac{\o^n_0}{n!} \leq C\Bigg(\frac{2s}{2s + \sigma}\Bigg)^2
	\end{equation}
	
	by \ref{RescaledKeyInequality}.
	
	Hence, we have
	\begin{equation}
		\int_X |(\Id - \H^s)\na \H^\sigma|^2 \frac{\o^n_0}{n!} \leq C\Bigg(\frac{2s}{2s + \sigma}\Bigg)^2,
	\end{equation} 
	
	which implies, after taking the limit $s \to 0$,
	\begin{equation}
		\int_X |\pi\na \H^\sigma|^2 \frac{\o^n_0}{n!} \leq 0
	\end{equation}

	and
	
	\begin{equation}
		\int_X |\pi\na (\Id - \pi)|^2 \frac{\o^n_0}{n!} \leq 0.
	\end{equation}
	
	This is because $\na \H^\sigma$ tends weakly to $\na (\Id - \pi)$. The fact that $\pi$ defines a coherent subsheaf follows from the separate meromorphicity theorem (see \cite{Shiffman, Popovici, UhlenbeckYau86, UY2}). The fact that $\mathcal{F}$ is nontrivial follows from the fact that $\det \tilde{h}_j$ tends to zero, and $\mathcal{F} \neq E$ follows from \ref{MoserEstimate}.
\end{proof}

\hfill

Finally we need a version of the theorem of Uhlenbeck and Yau, which shows us that our $\pi$ defines a subsheaf $\mathcal{F}$ of $E$, which destabilizes $E$ (with respect to $\o_0$). The difficulty comes from the fact that our metric is not fixed. Moreover, it varies with the metric $H_j$ in the fiber, so we cannot take the limit directly, as Uhlenbeck and Yau did. 

Let $\o_j$ be the metric, corresponding to $h_j$, and $V_j = \Vol(X, \o_j)$. We start with the standard Chern-Weil formula:
\begin{equation}\label{SlopeFormula}
	\mu_{\o_j}(\mathcal{F}) = {1 \over \text{rank}(\mathcal{F}) V_j}\Biggl\{\int_{X\setminus Z}  \Tr ((\sqrt{-1}\Lambda_{\o_j}F_{H_1} - \mu_{\o_j}(E) \Id)\pi) {\o^n_j \over n!} - \int_{X\setminus Z} |\na \pi|^2_{\o_j} {\o^n_j \over n!}\Biggr\}  + \mu_{\o_j}(E),
\end{equation}

and $Z=Z_\mathcal{F}$ is an analytic subset with $\codim_\C Z \geq 2$, where $\mathcal{F} \hookrightarrow E$ fails to be a subbundle.

We have for all $j$:
\begin{align*}
	&\Tr ((\sqrt{-1}\Lambda_{\o_j}F_{H_1} - \mu_{\o_j}(E) \Id)\pi) = \Tr ((\sqrt{-1}\Lambda_{\o_j}F_{H_j} - \mu_{\o_j}(E) \Id)\pi) + \Tr(g^{a \bar b}_j \p_{\bar b} (h^{-1}_j \na_a h_j)\pi) 
	\\
	&=\Tr(g^{a \bar b}_j \p_{\bar b} (h^{-1}_j \na_a h_j)(\pi - \Id)) + \Tr(g^{a \bar b}_j \p_{\bar b} (h^{-1}_j \na_a h_j)) 
	\\
	&= \Tr(g^{a \bar b}_j \p_{\bar b} (h^{-1}_j \na_a h_j)(\pi - \Id)) + g^{a \bar b}_j \p_{\bar b}\Tr( (h^{-1}_j \na_a h_j))
	\\
	&= \Tr(g^{a \bar b}_j \p_{\bar b} (h^{-1}_j \na_a h_j)(\pi - \Id)) + \Delta_j \log \det h_j.
\end{align*}

The second summand has zero integral over $X\setminus Z$.

So, we have the following formula:
\begin{equation}
	\mu_{\o_j}(\mathcal{F}) = {1 \over \text{rank}(\mathcal{F})V_j}\int_{X\setminus Z}\Biggl\{-\Tr(g^{a \bar b}_j \p_{\bar b} (h^{-1}_j \na_a h_j)(\Id - \pi)) - |\na \pi|^2_{\o_j}\Biggr\}{\o^n_j \over n!} + \mu_{\o_j}(E).
\end{equation}

Thus, we need to show that the following proposition holds:

\hfill

\begin{proposition}\label{SketchyProposition}
	The following inequality holds:
	\[\liminf_{j \to +\infty}\int_X\{\dots\}  \geq 0.\]
\end{proposition}

\hfill 

\begin{proof}
	
	We start with a simple rearrangement:	
\begin{align}
	&-\int_{X\setminus Z}\Tr(g^{a \bar b}_j \p_{\bar b} (h^{-1}_j \na_a h_j)(\Id - \pi)){\o^n_j \over n!} = -\int_{X\setminus Z}\Tr(g^{a \bar b}_j \p_{\bar b} (h^{-1}_j \na_a h_j)(\tilde{h}^\sigma_j)){\o^n_j \over n!}
	\\
	& + \int_{X\setminus Z} \Tr(g^{a \bar b}_j \p_{\bar b} (h^{-1}_j \na_a h_j)(\tilde{h}^\sigma_j - \Id + \pi)){\o^n_j \over n!} = \textbf{(I)} + \textbf{(II)}.
\end{align}

Let us deal with \textbf{(I)}. First, note that 
\[-\int_{X\setminus Z}\Tr(g^{a \bar b}_j \p_{\bar b} (h^{-1}_j \na_a h_j)(\tilde{h}^\sigma_j)){\o^n_j \over n!} \geq \int_{X\setminus Z} |\tilde{h}_j^{-\sigma/2}\na \tilde{h}_j^\sigma|^2{\o^n_j \over n!} \geq \int_{X\setminus Z} |\na \tilde{h}_j^\sigma|^2{\o^n_j \over n!} \]

by the proof of \ref{SobolevEstimate}.

\hfill

\begin{lemma}\label{SketchyLemma}
	We have the following inequality:
	
	\begin{equation}
		\liminf_{\sigma \to 0} \liminf_{j \to +\infty} \int_{X\setminus Z} |\na \tilde{h}^\sigma_j|^2_{\o_j}{\o^n_j \over n!} \geq \int_{X\setminus Z} |\na \pi|^2_{\o_0}{\o^n_0 \over n!}.
	\end{equation}
\end{lemma}

\hfill

\begin{proof}(of \ref{SketchyLemma})
	
	This is the place where we need the uniform convergence of $\o_s$ to $\o_0$ as $s \to 0$. The uniform convergence guarantees, that there are positive $\varepsilon_j \to 0$, such that
	\begin{equation}\label{SketchyLimits}
		(1 - \varepsilon_j)\int_{X\setminus Z} |\na \tilde{h}^\sigma_j|^2_{\o_0}{\o^n_0 \over n!} \leq \int_{X\setminus Z} |\na \tilde{h}^\sigma_j|^2_{\o_j}{\o^n_j \over n!} \leq (1 + \varepsilon_j)\int_{X\setminus Z} |\na \tilde{h}^\sigma_j|^2_{\o_0}{\o^n_0 \over n!}
	\end{equation}
	
	We immediately have 	
	\begin{align}
		\int_{X\setminus Z} |\na \H^\sigma|^2_{\o_0} {\o^n_0 \over n!} \leq \liminf_{j \to +\infty}\int_{X\setminus Z} |\na \tilde{h}^\sigma_j|^2_{\o_j}	{\o^n_j \over n!}
	\end{align}
	
	
	
	Taking the limit $\lim_{\sigma \to 0}$, we obtain the desired result.
\end{proof}

\hfill

Now we estimate \textbf{(II)}:
\begin{align}
	&|\textbf{(II)}| \leq \int_{X\setminus Z}|\Tr (g^{a\bar b}_j F_{a \bar b}(\tilde{h}^\sigma_j - \Id + \pi))|{\o^n_j \over n!} + \int_{X\setminus Z}|\Tr (g^{a \bar b}_j F_{1a \bar b}(\tilde{h}^\sigma_j - \Id + \pi))|{\o^n_j \over n!}
	\\
	& \leq C \int_{X\setminus Z}|\tilde{h}^\sigma_j - \Id + \pi| {\o^n_0 \over n!},
\end{align}

and the RHS tends to zero, because $\tilde{h}^\sigma_j$ tends a.e. to  $\H^\sigma$ as $j \to +\infty$, and $\H^\sigma$ tends a.e. to $\Id - \pi$, as $\sigma \to 0$.

This concludes the proof of \ref{SketchyProposition}.
\end{proof}

We have shown that 
\[\lim_{j \to +\infty}\mu_{\o_j}(\mathcal{F}) \geq \mu_{\o_0}(E).\]

In fact, in our setting, we also have 
\[\lim_{j \to +\infty}\mu_{\o_j}(\mathcal{F}) \leq \mu_{\o_0}(E),\]
and this inequality holds for any subsheaf of $E$, since  for any $s>0$ our bundle $E$ is $\o_s$-stable. Hence we have
\[\mu_{\o_0}(\mathcal{F}) = \mu_{\o_0}(E).\]
Thus, combining these inequalities, we obtain the following proposition, which is the second part of \ref{MainTheorem}.

\hfill

\begin{proposition}\label{MainSecondPart}
	If the family of relative endomorphisms $h_s$, normalized by 
	\[\int_X \log \det h_s \o_s^n = 0,\] 
	has unbounded traces $\Tr h_s$, then the bundle $E$ is strictly $\o_0$-semistable. Moreover, a projection onto a destabilizing sheaf may be constructed by the Uhlenbeck-Yau construction. 
\end{proposition}

\hfill

\begin{remark}
	Our proof carries over verbatim even if the family of K\"ahler metrics $\o_s$ in our considerations is replaced by a family of balanced (i.e. $d \o_s^{n-1} = 0$) or Gauduchon (i.e. $\p \bar{\partial}\o_s^{n-1} = 0$) metrics. The only place that should be modified is the proof of \ref{NormalizationOfEnd}. Since the uniform equivalence allows us to control the Sobolev constant, the usual Moser iteration argument shows that $\log \det h_s$ is uniformly bounded. 
\end{remark}

\hfill

\begin{remark}
	Our arguments show, that the metrics $\tilde{H}_j = {H_j \over \sup_X \Tr h_j}$ degenerate to some limiting metric on $E / \mathcal{F}$ over $X \setminus Z$. However, the convergence of these metrics is very weak, as one can see directly from our arguments above. In particular, $H_j$ and $\tilde{H}_j$ have the same Chern connections $A_{H_j}$, and curvatures $F_{H_j}$.
	We address the behaviour of the connections in the next section.  
\end{remark}

\section{Uhlenbeck limits of connections}\label{UhlenbeckCompSection}

In this section we discuss the Uhlenbeck limits of the Chern connections of metrics $H_s$.
While we have shown that the metrics themselves blow up, their connections still may have a limit (modulo gauge). However, there is a delicate point here: while there is a well-developed compactness theory for the Yang-Mills connections on Riemannian manifolds, this theory works under the assumption that the background K\"ahler (or just a Riemannian) metric is fixed. Our setting is different, as we allow the base metrics $\o_s$ to vary, so we shall show that Uhlenbeck compactness and other results hold in the setting of varying metrics as well. In order to do it, we need to impose the following conditions on $\o_s$:
\begin{enumerate}
	\item $\o_s$ are uniformly equivalent to $\o_0$;
	
	\item $\o_s$ converge in $C^0$ to $\o_0$, as $s \to 0$;
	
	\item There is a constant $C>0$, such that $|\text{Rm}(\o_s)|_{\o_s} \leq C$.
\end{enumerate}

The only new condition here is the bound for curvature. At some points, the bound for the curvature enters our estimates. We also need to control the uniform lower bound for the injectivity radii of $\o_s$, and the bound for the curvature together with the first condition (which translates to uniform volume and diameter bounds) allow us to obtain such a bound.

Also, for higher regularity estimates, we need to control $|\na^k\text{Rm}(\o_s)|$. We will always state explicitly, when such bounds are required for our purposes.

\subsection{Preliminaries: $\varepsilon$-regularity, monotonicity formula and the estimates of curvature}

In this section we turn our attention to the limiting behaviour of the Chern connections of our Hermitian-Einstein metrics. First of all, we need to recall the monotonicity formula for the Yang-Mills connections, which we slightly reformulate for our purposes.

\hfill

\begin{theorem}\label{MonotonicityFormulae}(see \cite[Thm. 2.1.2]{Tian2000} and \cite{Price83})
	Let $(X, \o_s)$ be a K\"ahler manifold with a family of K\"ahler metrics $\o_s$. Assume that the metrics are equivalent, and their curvatures are uniformly bounded. Let $(E,H)$ be a hermitian vector bundle over $X$, and $A_s$ be a family of the integrable Hermitian Yang-Mills connections with respect to $\o_s$. Then there is a positive number $R_0$ and $a \geq 0$, such that for all $0 < r < R \leq R_0$ and for any $s$ the following holds:
	\begin{align}
		&R^{4-2n}e^{aR^2}\int_{B_{R,s}(p)} |F_{A_s}|^2 {\o^n_s \over n!} - r^{4-2n}e^{ar^2}\int_{B_{r,s}(p)} |F_{A_s}|^2 {\o^n_s \over n!} 
		\\
		& \geq 2 \int_{B_{R,s}(p) \setminus B_{r,s}(p)} \rho_s^{4-2n}e^{a\rho^2_s} |\iota_{\p_{\rho_s}}F_{A_s}|^2{\o^n_s \over n!}
	\end{align}
	where $B_{r,s}(p)$ denotes the ball of radius $r$ in $(X, \o_s)$, and $\rho_s$ is the distance function with respect to $\o_s$. Moreover, the equality holds if $X = \C^n$ with standard flat metric and $a=0$.
\end{theorem}

\hfill

Next, we need the $\epsilon$-regularity statement.

\hfill

\begin{proposition}\label{EpsilonReg}
	Let $(X, \o_s, E, H_s)$ be as in our setting. There are $\epsilon_0 >0, ~C>0$ and $R_0>0$, such that if for $\epsilon \leq \epsilon_0$ and $0<r<R_0$ we have 
	
	\[r^{4-2n}\int_{B_{r,s}(p)} |F_{H_s}|^2 {\o^n_0\over n!} < \epsilon, \]
	
	then 
	
	\[|F_{H_s}|(p) \leq {C \over r^2} \Bigg(r^{4-2n}\int_{B_{r,0}(p)} |F_{H_s}|^2 {\o^n_0 \over n!}\Bigg)^{1 \over 2}.\]
\end{proposition}

\hfill

\begin{proof}
	Since $F_{H_s}$ satisfies the (Hermitian-)Yang-Mills equations (with respect to $\o_s$), and we have the uniform control of $|\text{Rm}(\o_s)|$, we have
	
	\[\Delta_{\o_s} |F_{H_s}| \geq - C_1 (|F_{H_s}|^2 + |F_{H_s}|),\]
	
	for some universal $C_1 >0$ independent of $s$. Also, the monotonicity formula for the Yang-Mills connections holds (for fixed $s$), so does the $\varepsilon$-regularity. Since $\o_s$ is uniformly equivalent to $\o_0$, we can replace $\o_s$ by $\o_0$ by adjusting the constant on the right hand side. 
\end{proof}

\hfill

\subsection{Uhlenbeck compactness.}

In this section we investigate the limits of Chern connections $\na^s$ of the Hermitian-Einstein metrics $H_s$, as $s$ approaches $0$. First of all, we note that the Uhlenbeck compactness theorem does not work in our setting, since each $\na^{H_s}$ is hermitian, but with respect to its own metric $H_s$. To fix this problem, we introduce $w_s = h^{1/2}_s$, and introduce new connections
\[\na_{A_s} = d + A_s = w_s \circ \na^{H_s} \circ w^{-1}_s.\]
These connections are integrable and unitary with respect to fixed background metric $H = H_1$. Moreover, each $A_s$ is Hermitian Yang-Mills with respect to $\o_s$. These connections enjoy the following obvious properties.

\hfill

\begin{lemma}\label{GaugeIsometric}
	We have the following formulas
	\begin{align*}
		&|F_{H_s}|^2_{H_s} = |F_{A_s}|^2_H
		\\
		&|\Lambda_\o F_{H_s}|^2_{H_s} = |\Lambda_\o F_{A_s}|^2_H,
	\end{align*}
	and the same holds for covariant derivatives of $F_{A_s}$. Here $\o$ is any fixed K\"ahler form on $X$.
\end{lemma}

\hfill

The proof follows from the fact that $w_s$ is a bundle isometry from $(E, H_s)$ to $(E, H)$, so we omit it. In the following we drop subscript $H$ from $|F_{A_s}|^2_H$ to lighten the notation.

However, another difficulty arises in our course. In our setting, we have uniform control of $|\Lambda_{\o_s} F_{A_s}|$, but we do not control $|\Lambda_{\o_0} F_{A_s}|$. Even if we have fairly good control of our family of metrics, we cannot directly convert this into the estimate for $|\Lambda_{\o_0} F_{A_s}|$ in terms of $\o_s, \o_0$ and $|\Lambda_{\o_s} F_{A_s}|$. So we cannot apply the Uhlenbeck compactness theorem directly in the form as it is stated, for example, in \cite[Theorem 2.17]{Sibley}. So, we need to show, that for the family $\{A_s\}$ there is a version of Uhlenbeck compactness theorem. 

Before we proceed, we need the following definition.

\hfill

\begin{definition}\label{ZanDefinition} Let $\{s_j\}$ be a sequence, such that $s_j \to 0$ as $j \to +\infty$. Write $A_j = A_{s_j}$, and the corresponding K\"ahler metric $\o_j = \o_{s_j}$. 
	\begin{equation}
		Z_{an} = \Biggl\{p \in X ~\vline~  \lim_{r \to 0}\liminf_{j \to +\infty} r^{4-2n}\int_{B_{r,j}(p)} |F_{A_j}|^2 {\o^n_j \over n!} \geq \epsilon_0 ~\Biggr\}
	\end{equation}
	
	Here $B_{r,j}(p)$ is a ball in $(X, \o_j)$ of radius $r$. 
\end{definition}

\hfill

Of course, a priori, $Z_{an}$ depends on the chosen sequence. However, we omit the indication of such dependence, because it will be clear, which sequence we are working with.

We have the following proposition.

\hfill

\begin{proposition}\label{PropertiesOfZan}
	One can replace $\o_s$ in this definition with $\o_0$. Moreover, $Z_{an}$ is closed and has Hausdorff dimension at most $2n-4$. 
\end{proposition}

\hfill

\begin{proof}
	The statement follows from our assumption, that any sequence $\o_j$ converges to $\o_0$ uniformly, and metrics are uniformly equivalent. In particular, $B_{r,j}(p)$ converges to $B_{r,0}(p)$ for any $p$. In particular, there are sequences of positive constants $\{\ve_j\},\{\delta_j\}$, both tending to zero as $j \to +\infty$, such that
	\[(1-\ve_j)\int_{B_{r-\delta_j,j}(p)}|F_{A_j}|^2 {\o^n_0 \over n!} \leq \int_{B_{r,0}(p)}|F_{A_j}|^2 {\o^n_0 \over n!}\leq (1+\ve_j) \int_{B_{r + \delta_j,j}(p)}|F_{A_j}|^2 {\o^n_j \over n!}.\]
	
	The statement about the Hausdorff dimension of $Z_{an}$ is standard, see for instance \cite{Nakajima88, Tian2000}.  Q.E.D. 
\end{proof}

\hfill

Now we can state the version of Uhlenbeck compactness theorem, which works in our setting. 

\hfill

\begin{theorem}\label{UhlenbeckCompactness}
	Let $X$ be a compact complex manifold, and let $\{\o_s\}_{s \in [0;1]}$ be a family of K\"ahler metrics, which are uniformly equivalent to $\o_0$, converge uniformly to $\o_0$, and have uniformly bounded curvature. Let $(E,H)$ be a hermitian vector bundle, and suppose $\{A_s\}$ is a family of integrable hermitian connections. Assume that $||F_{A_s}||_{L^2(X, \o_0)} $ and $|\Lambda_{\o_s}F_{A_s}|$ are uniformly bounded.  Fix $p >2n$. Then, there is a sequence $s_j \to 0$, such that the following holds for $A_j = A_{s_j}$:
	\begin{enumerate}
		\item There is a $C^1$ hermitian vector bundle $(E_\infty, H_\infty)$ on $X \setminus Z_{an}$, where $Z_{an}$ is defined in \ref{ZanDefinition};
		
		\item A hermitian connection $A_\infty$ on $E_\infty$ with finite $\YM (A_\infty)$, such that
		$A_j$ converge to $A_\infty$ modulo gauge transforms weakly in $L^p_{1, loc}(X \setminus Z_{an})$.  
		
	\end{enumerate}
\end{theorem}

\hfill

\begin{proof}
	As shown in \cite[pages 4 and 13]{UhlenbeckPreprint}, the metric $\o$ enters in estimates for connection by the lower bound to its injectivity radius from below, and the bound for its curvature. By our assumptions on metrics $\o_s$, we have the uniform control of their curvatures. Also, by uniform equivalence, we control volumes and diameters of $\o_s$. By Cheeger Finiteness Theorem, we have uniform lower bound on the injectivity radii of $\o_s$, hence estimates from \cite[Chapter 2]{UhlenbeckPreprint} are uniform in $s$. Now we apply arguments of Uhlenbeck and the standard elliptic theory and obtain the sequence of $A_j$ with stated properties. 
\end{proof}

\hfill

In our setting, where $A_s$ are Hermitian Yang-Mills connections constructed from the Chern connections of Hermitian-Einstein metrics, we can say more.

\hfill

\begin{corollary}\label{UhlenbeckLimitOfYM}
	If in addition $\sqrt{-1}\Lambda_{\o_s}F_{A_s} = \mu_{\o_s}(E)\Id$, then the Uhlenbeck limit $A_\infty$ is Hermitian Yang-Mills with respect to $\o_0$. Moreover, $A_j$ converge to $A_\infty$ in $C^\infty_{loc}(X \setminus Z_{an})$, if for every $k$ we have the uniform $C^k$-bounds for $\text{Rm}(\o_s)$, and $E_\infty$ extends over all of $X$ as a reflexive sheaf.
\end{corollary}

\hfill

\begin{proof}
	The fact that the limit is Yang-Mills follows from taking the weak limit of connections (modulo gauge) in $L^p_1, ~p>2n$. Indeed, the assumption $\sqrt{-1}\Lambda_{\o_j}F_{A_j} = \mu_{\o_j}(E)\Id$ implies that
	\[\sqrt{-1}\Lambda_{\o_0}F_{A_\infty} = \mu_{\o_0}(E)\Id\]
	in the sense of distributions. 
	The statement about the $C^\infty$ convergence can be shown along the lines of \cite{Nakajima88}, and the fact that $A_\infty$ is Hermitian Yang-Mills follows from it directly.
	The fact that $E_\infty$ extends to a reflexive sheaf over all of $X$ follows from \cite{BandoSiu}.
\end{proof}

\hfill

However, we can say more about $A_\infty$. Now we show the main result of this section.

\hfill

\begin{theorem}\label{ReductionToTheYMFlow}
	The Uhlenbeck limit $A_\infty$ defines a reflexive sheaf $E_\infty$, which is isomorphic to $\text{Gr}^S_{\o_0}(E)^{**}$.   
\end{theorem}

\hfill

\begin{proof}
	Consider the following YM flow on $(X, \o_0)$:
	\begin{equation}
		\left\{\begin{array}{@{}l@{}}
			\p_t A_{t,j} = - d^*_{A_{t,j}}F_{A_{t,j}}
			\\
			A_{t,j}|_{t=0} = A_j.
		\end{array}\right.\,.
	\end{equation}
	Since $\lim_{j \to + \infty}\YM_{\o_j}(A_j) = \lim_{j \to + \infty} \YM_{\o_0}(A_j)$, we know that $\YM_{\o_0}(A_j)$ tends to the infimum of $\YM_{\o_0}$ by \ref{InfAttained}. Now the same argument, as in \cite{DasWentworth1} and \cite[Proposition 7.2]{Sibley} shows that
	\begin{align*}
		&||A_{t,j} - A_j||_{L^2(X, \o_0)} \leq \int_{0}^{t}||\p_s A_{s,j}||_{L^2(X, \o_0)}ds 
		\\
		&\leq \sqrt{t}\Bigg(\int_{0}^{t}||d^*_{A_{s,j}}F_{A_{s,j}}||_{L^2(X, \o_0)}^2ds\Bigg)^{1/2}
		\\
		&\leq \sqrt{{t \over 2}}(\YM_{\o_0}(A_j) - \YM_{\o_0}(A_{t,j}))^{1/2} \xrightarrow[j \to +\infty]{} 0,
	\end{align*}
	because $\YM_{\o_0}(A_j) \to \inf_{A \in \mathcal{A}^{1,1}_H}\YM_{\o_0}(A)$ as $j \to + \infty$ by \ref{InfAttained}, and $\YM_{\o_0}$ is decreasing along the Yang-Mills flow.
	In particular, after taking a further subsequence, we have
	\[||d^*_{A_{s,j}}F_{A_{s,j}}||_{L^2(X, \o_0)} \xrightarrow[j \to + \infty]{} 0\]
	for almost all $s$. Hence we showed, that the Uhlenbeck limit of $A_{t_0,j}$ is the same as the limit of $A_j$, and we again showed that the limit of $A_j$ is Hermitian Yang-Mills with respect to $\o_0$. Also, the standard argument shows that
	\[|\Lambda_{\o_0} F_{t,j}|^2(x) \leq \int_X K_t(x, y) |\Lambda_{\o_0} F_{A_j}|^2 {\o^n_0 \over n!},\]
	where $K_t(x,y)$ is the heat kernel of the Laplacian on $X$ with respect to $\o_0$.
	Since $K_t(x,y) \leq C(1 + t^{-n})$ for any $t>0$, the function $|\Lambda_{\o_0} F_{t,j}|^2$ is uniformly bounded above by a constant times $||\Lambda_{\o_0}F_{A_j}||^2_{L^2(X,\o_0)}$. 
	Also, note that for any $j$ the operators $\bar{\partial}_{A_j}$ lie in the same complex gauge orbit. Hence, the operators $\bar{\partial}_{A_{t_0,j}}$ also lie in the same gauge orbit.
	
	Finally, from
	\[\YM_{\o_0}(A) = \HYM_{\o_0}(A) + {4\pi^2 \over (n-2)!}\int_X(2c_2(E) -c^2_1(E)) \wedge\o^{n-2}_0,\]
	we see that 
	\[\lim_{j \to + \infty}\HYM_{\o_j}(A_{t,j}) = \lim_{j \to + \infty}\HYM_{\o_0}(A_j) = \HYM_{\o_0}(\mu_0),\]
	where $\mu_0 = (\mu_0(E), \dots, \mu_0(E))$ -- the Seshadri type of $E$. Hence, the same is true for $\HYM_{\alpha, N,\o_0}(A_j)$ for any $\alpha \in[1;2]$, and any $N>0$. In other words, for any $\alpha \in [1;2]$ we have
	\[{1 \over r}\sum_{\gamma = 1}^{r} |\mu_\gamma(A_\infty) + N|^\alpha = |\mu_0(E) + N|^\alpha.\]
	Hence for any $\gamma$ we obtain $\mu_\gamma(A_\infty) = \mu_{\o_0}(E)$. 
	
	We have shown that the limit of $A_j$ is the same as the limit of $A_{t_0,j}$ for some $t_0$. Moreover, the Harder-Narasimhan type of $(E_\infty, A_\infty)$ is the same as $\mu_0$. Now the arguments from \cite{Donaldson85}, or \cite[Section 7.1]{Sibley}, or invoking \cite[Theorem 2.7]{SibleyWentworth} show the isomorphism of $E_\infty$ to $\text{Gr}^S_{\o_0}(E)^{**}$. 
\end{proof}

\subsection{Rectifiability and analyticity of the blow-up set.}

At this moment, we do not know that $Z_{an}$ is indeed an analytic subvariety of $X$. We do not even know whether it is rectifiable. So, in this part we are going to establish these properties of $Z_{an}$. For the rectifiability part, we will mimic the exposition in \cite[Chapter 3]{ChenWentworth1}, and for analyticity we will follow the arguments of Hong and Tian from \cite{HongTian} and \cite{Tian2000}. First of all, we need the following definition

\hfill

\begin{definition}\label{DefectMeasure}	
	Let $\{s_j\}$ be a sequence, which tends to zero, and $\o_j$ with $A_j$ are the corresponding sequences of metrics and unitary connections. Then one can define a Radon measure $\mu$ as
	
	\[\mu = \lim_{j \to +\infty}|F_{A_j}|^2 {\o^n_{j} \over n!} = |F_{A_\infty}|^2 {\o^n_0 \over n!} + \nu.\]
	
	The measure $\nu$ is called the \textbf{defect measure}.
\end{definition}

\hfill

Of course, this measure depends on the sequence $A_j$.

We need the following proposition

\hfill

\begin{proposition}\label{LongProposition} 
	
	Let $(X, \o_{j}, E, A_j)$ be as above. Let $B_r(x)$ be the ball of radius $r$ with respect to $\o_0$. Then the following holds
	\begin{enumerate}
		\item For $0<r<<1$ the quantity $e^{ar^2}r^{4-2n}\mu(B_{r}(x))$ is decreasing as $r \to 0$. Hence, the function 
		\[\Theta(\mu,x) = \lim_{r \to 0}r^{4-2n}e^{ar^2}\mu(B_{r}(x))\]
		is well-defined, upper semi-continuous and $\mathcal{H}^{2n-4}$ approximately continuous at $\mathcal{H}^{2n-4}$ a.e. $x \in Z_{an}$;
		
		\item  The point $x$ belongs to $Z_{an}$ if and only if $\Theta(\mu,x) \geq \epsilon_0$.  Moreover, for $\mathcal{H}^{2n-4}$ a.e. $x \in Z_{an}$ we have that
		\[\limsup_{r \to 0} r^{4-2n}e^{ar^2}\int_{B_{r}(x)}|F_{A_\infty}|^2 \o^n_{0} = 0.\]
		
		\item The defect measure $\nu$ is $\mathcal{H}^{2n-4}_{Z_{an}}$ absolutely continuous and is equal to $\Theta(\mu, x)d\mathcal{H}^{2n-4}_{Z_{an}}$.
	\end{enumerate}
	
\end{proposition}

\hfill

\begin{proof}
	The arguments from \cite{ChenWentworth1, HongTian, FLin, FLinBook, Tian2000} can be carried over almost verbatim here. We just outline the proof, and explain the main differences. 
	
	\begin{enumerate}
		\item This follows from the observation, that $|F_{A_j}|^2 \o^n_{j}$ and $|F_{A_j}|^2 \o^n_{0}$ have the same weak limit $n!\mu$. Indeed, our assumptions on $\o_s$ imply that $\o_j \to \o_0$ uniformly. Hence, there are $\delta_j \to 0$, such that:
		
		\[\lim_{j \to +\infty}\int_{B_{r + \delta_j,j}(x)}|F_{A_j}|^2\o^n_{j} = \lim_{j \to +\infty}\int_{B_{r}(x)}|F_{A_j}|^2 \o^n_{0},	\]
		
		again because $\o_j$ are uniformly equivalent and converge to $\o_0$ uniformly. By \ref{MonotonicityFormulae} we see that $r^{4-2n}e^{ar^2}\mu(B_{r,\o_0}(x))$ is non-decreasing for every $x$, and the limit exists. We only need to ensure, that $\mu(\p B_{r,0}(x)) = 0$, as measure may concentrate on the spheres as $j \to + \infty$. However, this problem does not exist for almost all radii. 
		
		The upper-semicontinuity is a standard fact, and the approximate continuity follows from the argument of \cite[Lemma 3.2.2]{Tian2000} almost verbatim. 
		
		\item The statement follows from standard covering arguments, see \cite{Tian2000} or \cite{FLinBook}. 
		
		\item The proof is similar to \cite[Lemma 3.3]{ChenWentworth1}. 
	\end{enumerate}
\end{proof}

\hfill

We also need the following proposition, which is an analogue of decomposition from \cite[Lemma 3.3]{ChenWentworth1}.

\hfill

\begin{proposition}\label{StructureOfBlowUp1}
	The set $Z_{an}$ can be decomposed as 
	\[Z_{an} =  Z_{b, an} \cup \text{Sing}(A_\infty),\]
	where $Z_{b, an} = \text{supp}(\nu)$, and 
	\[\text{Sing}(A_\infty) =\{x \in X ~|~ \limsup_{r \to 0} r^{4-2n}\int_{B_{r,0}(x)} |F_{A_\infty}|^2 \o^n_0 \geq \ve_0\}\] 
	is an analytic subset of codimension at least 3. 
\end{proposition}

\hfill

\begin{proof}
	This follows essentially from (3) of \ref{LongProposition}, and it can be proved along the lines of \cite[Lemma 3.3]{ChenWentworth1}, by using the fact that metrics $\o_j$ are uniformly equivalent to $\o_0$. The fact that $\text{Sing}(A_\infty)$ is analytic and of codimension at least 3 follows from the result of Bando and Siu \cite{BandoSiu}.
\end{proof}

\hfill

It can be shown, that one can find suitable local gauge transformations $u$, such that $u^*A_{\infty}$ extends across $Z_{b,an} \setminus \text{Sing}(A_\infty)$. 

Next, we need to consider tangent cone measures. Before we proceed, we need to set some notation.

Let $\tau_\lambda: T_x X \to T_x X$ be a dilation by $\lambda$. Let $\mu_{x,\lambda}(E) = \lambda^{4 -2n}\mu(\exp_{x,g_0}(\tau_\lambda (E)))$ for any $E \subset T_x X$. Here $\exp_{x,g_0}$ is the exponential map of the metric $\o_0$. Denote also by $g_j$ the Riemannian metric for $\o_{s_j}$, and by $g_0$ the Riemannian metric given by $\o_0$. We have the following proposition

\hfill

\begin{proposition}\label{TangentConeExists}
	There is a sequence $\lambda_j \to 0^+$, such that $\mu_{x,\lambda_j} \to \eta$, where $\eta$ is a cone measure on $T_x X$, i.e.
	
	\[\eta(\lambda E) = \lambda^{2n-4}\eta(E).\]  
\end{proposition}

\hfill

\begin{proof} 
	
	This essentially follows from \ref{LongProposition}, together with Preiss' theorem from \cite{Preiss}. However, we can give a short proof, following \cite[Proposition 3.4]{ChenWentworth1}.
	
	Consider the following connections:
	\[A_{x,j,\lambda} = \tau^*_{\lambda}\exp^*_{x, g_0}A_{s_j}.\]
	
	Let $g_{x,j,\lambda} = {1 \over \lambda^2}\tau^*_\lambda\exp^*_{x, g_0}g_j$, and let $dV_{x,j,\lambda}$ be the volume form with respect to this metric.
	
	Let $\eta = \lim_{\lambda_k \to 0}\mu_{x,\lambda_k}$. This limit exists because of the monotonicity formula from \ref{MonotonicityFormulae}. There is a subsequence $A_{x,j,\lambda_j}$ of $A_{x,j,\lambda_k}$, such that
	\[|F_{x,j,\lambda_j}|^2 dV_{x,j,\lambda_j} \xrightarrow[j \to +\infty]{} \eta.\]
	
	Note that $g_{x,j,\lambda_j}$ and $g_{x,0, \lambda_j}$ have the same limit, as metrics are uniformly equivalent, and $g_j$ tends to $g_0$ as $j \to +\infty$. The fact that $\eta$ is a cone measure now follows from the monotonicity formula \ref{MonotonicityFormulae}, just as in \cite[Proposition 3.4]{ChenWentworth1} or \cite[Lemma 3.2.1]{Tian2000}, but with a little modification. For fixed $j$, and any radially invariant function $\phi \geq 0$, we have
	\begin{align}\label{ConeFormula}
		&R^{4-2n}\int_{B_{R,g_j}(x)} \phi |F_{x,j,\lambda_j}|^2 dV_{x,j, \lambda_j} - r^{4-2n}\int_{B_{r,g_j}(x)} \phi |F_{x,j,\lambda_j}|^2 dV_{x,j,\lambda_j} 
		=
		\\
		&\int_{r}^{R}dt {d \over dt}\Bigg(t^{4-2n}\int_{B_{1,g_j}(x)}\phi |\tau^*_t F_{x,j,\lambda_j}|^2 \tau^*_tdV_{x,j,\lambda_j}\Bigg)
	\end{align}
	
	 As in \cite[Proposition 3.4]{ChenWentworth1}, the expression on the right hand side of \ref{ConeFormula} is bounded by a constant times the integral of
	\[r^{4-2n}|\iota_{\p_r} F_{x,j,\lambda_j}||F_{x,j,\lambda_j}|.\]
	
	Here $r=r_j$ -- the distance with respect to $g_j$. As $\lambda_j \to 0$, we know that $g_{x,j,\lambda_j}$ tend uniformly to the Euclidean metric $\delta = g_{x,0}$, so the integral of $|\iota_{\p_r} F_{x,j,\lambda_j}|$ tends to zero by \ref{MonotonicityFormulae}. 
	
\end{proof}

\hfill

Next we show the existence of weak complex tangent planes.

\hfill

\begin{proposition}\label{TangentPlaneExists}
	The set $Z_{an}$ is rectifiable, and for $\mathcal{H}^{2n-4}$ a. e. $x \in Z_{an}$ the  tangent plane to $Z_{an}$ at  $x$ is a complex subspace in $T_x X$, calibrated by ${\o^{n-2}_0 \over (n-2)!}$.  
\end{proposition}

\hfill

\begin{proof}
	The rectifiability follows from \ref{LongProposition} together with the main result of \cite{Preiss}.
	For the last statement, that tangent planes are complex, the arguments from \cite[Section 4]{Tian2000} can be carried over word by word in our setting.
	
	Let $\tau^{\delta_j}_{y_1,y_2}(x) = \delta_jx + (y_1,y_2)$. Here $y_1 \in V$, and $y_2 \in V^\perp$. The forms
	\[\Omega_j = \delta_j^{4-2n}\tau^{\delta_j,*}_{y_1,y_2}\exp^*_x {\o^{n-2}_j \over (n-2)!}\]
	converge to $\Omega_0 = {\o_0^{n-2} \over (n-2)!}\vline_{T_x X}$. Tian's argument shows that there exists $\Omega_0$-anti-self-dual instanton $B$ with $F_B \neq 0$, $\Lambda F_B = 0$, and $\iota_v F_B = 0$ for any $v \in V$. We now can repeat the argument of \cite[Corollary 4.2.2]{Tian2000}, that $\Omega_0 |_{T_x Z_{an}} $ is equal to the volume form on $T_x Z_{an}$, defined by the metric $g_{x,0}$. In other words, $T_x Z_{an}$ is $\Omega_0$-calibrated, and by the Wirtinger inequality, it is a complex subspace.
\end{proof}

\hfill

Now we are in a position to prove the analogue of the theorem of Hong and Tian.

\hfill

\begin{theorem}\label{ZanIsAnalytic}
	(see \cite[Theorem 11]{HongTian} and \cite[Theorem 4]{Tian2000}).
	
	Let $X$ be a compact K\"ahler manifold, such that the K\"ahler metrics $\o_j$ uniformly converge to $\o_0$, and have uniformly bounded curvature. Let $(E,H)$ be a hermitian vector bundle and suppose $A_j \in \mathcal{A}^{1,1}_H(E)$ are the integrable Hermitian Yang-Mills connections with respect to $\o_j$. Then, up to taking a subsequence, $A_j$ converge weakly to the $\o_0$ Hermitian Yang-Mills connection with the blow-up locus $Z_{an}$. The subset $Z_{an}$ is rectifiable and admits a decomposition $Z_{an} = \cup_\alpha Z_\alpha \cup \text{Sing}(A_\infty)$, where each $Z_\alpha$ is a codimension $2$ analytic subvariety. Moreover, for each smooth $(n-2,n-2)$ form $\phi$
	\[\lim_{j \to + \infty}\int_X c_2(A_j) \wedge \phi = \int_X c_2(A_\infty) \wedge \phi + \sum_\alpha m_\alpha \int_{Z_\alpha} \phi.\]
	
	Here for any connection $A \in \mathcal{A}^{1,1}_H$ we define $c_2(A)$ as 
	\[c_2(A) = {1 \over 8\pi^2}\Big(\Tr(F_A \wedge F_A) - \Tr F_A \wedge \Tr F_A\Big).\] 
\end{theorem}

\hfill

\begin{proof}
	The proof is essentially the same, as \cite[Theorems 4.2.3-4.3.3]{Tian2000}. We outline the main idea, and the modifications to add in our setting, so Tian's arguments would work.
	
	First of all, we need to show, that $\Tr F_{A_j} \wedge \Tr F_{A_j}$ converge to $\Tr F_{A_\infty} \wedge \Tr F_{A_\infty}$, and the latter extends to the whole $X$. Note that $\alpha_j = \sqrt{-1}\Tr F_{A_j}$ is $\o_j$-harmonic, and lies in $2\pi c_1(E)$. The Bochner formula and the Moser iteration imply that $C^0$-norm of $\alpha_j$ is controlled by its $L^2(X,\o_0)$-norm. Moreover, for any $\o$ the formula
	\[\alpha^2 \wedge {\o^{n-2} \over (n-2)!} = \Big(|\Lambda_\o\alpha|^2 - |\alpha|_\o^2\Big){\o^n \over n!},\]
	together with \ref{InfAttained} implies that the $L^2$-norm of $\alpha_j$ is uniformly bounded. Thus, by the standard elliptic theory, the $C^1$-norm of $\alpha_j$ is uniformly bounded by a constant, independent of $j$.
	
	Let $\alpha_0$ be a $\o_0$-harmonic form in $2\pi c_1(E)$. The same argument as above implies that
	\[||\alpha_j - \alpha_0||^2_{L^2(X, \o_0)} = ||\Lambda_{\o_0}\alpha_j - \Lambda_{\o_0}\alpha_0||^2_{L^2(X,\omega_0)} \xrightarrow[j \to + \infty]{}0.\]
	Clearly $\alpha_0 = \sqrt{-1}\Tr F_{A_\infty}$ on $X \setminus Z_{an}$. Hence, there is a bounded extension of $\Tr F_{A_\infty}$ on the whole $X$, and $\Tr F_{A_j} \wedge \Tr F_{A_j}$  converge to $\Tr F_{A_\infty} \wedge \Tr F_{A_\infty}$ as currents. 
	
	Consider the currents $T_j = c_2(A_j) - c_2(A_\infty)$. The previous discussion and the same arguments as \cite[Proposition 2.3.1]{Tian2000} imply, that $T_j$ are closed. The uniform $L^2$-bound on curvatures $F_{A_j}$ implies that $T_j$ have bounded mass. Thus, up to taking a subsequence, there is a limit current
	\[T = \lim_{j \to + \infty} T_j = \lim_{j \to + \infty} {1 \over 8\pi^2}\Big( \Tr (F_{A_j} \wedge F_{A_j})  - \Tr (F_{A_\infty} \wedge F_{A_\infty}) \Big), \]
	which is closed and has bounded mass. Up to replacing $g_{x,\lambda_j}$ by $g_{x,j,\lambda_j}$, and $\tau^*_{\lambda_j}\Omega$ by $\Omega_j$ from \ref{TangentPlaneExists}, we can repeat Tian's arguments, and show that $T$ is integer.
	
	Now, Simon's rectifiability theorem (see \cite[Theorem 32.1]{SimonGMT}) guarantees, that there is a rectifiable set $\Sigma$ of Hausdorff codimension 4, together with an orientation multivector $\eta$, and an integer-valued function $\theta$, such that
	\[T(\phi) = \int_\Sigma \langle \phi, \eta \rangle \theta d\mathcal{H}^{2n-4}_\Sigma.\]
		
	On the other hand, for $\phi = f {\o_0^{n-2} \over (n-2)!}$ we clearly have
	
	\[T(\phi) = {1 \over 8\pi^2}\int_{Z_{an}} f\Theta d\mathcal{H}^{2n-4}_{Z_{an}}.\]

 	Finally, as in Tian's paper, the results from \cite{King, HarveyShiffman} imply, that there are holomorphic subvarieties $Z_\alpha$, and positive integers $m_\alpha$, such that
	\[T(\vp) = \sum_\alpha m_\alpha \int_{Z_\alpha} \vp.\]
	
	Hence, $Z_{b, an} = \cup_\alpha Z_\alpha$, and $Z_{an} = Z_{b,an} \cup \text{Sing}(A_\infty)$ is analytic.
\end{proof}

\hfill

\begin{definition}\label{ANMultiplicities}
	We will refer to the numbers $m_\alpha$ that appear in the proof of \ref{ZanIsAnalytic} as {\em analytic multiplicities} of $Z_\alpha$.
	
\end{definition}

\subsection{The analog of Sibley-Wentworth theorem}

In this section we will prove the analogue of \cite[Theorem 1.1]{SibleyWentworth}. This result is the content of the last part of \ref{MainTheorem}. In fact, we just explain how the arguments of Sibley and Wentworth can be carried over in our framework.  

Before we proceed, we recall the notion of algebraic multiplicities. We need to introduce the notion of associated graded object of a coherent sheaf (see \cite[Erratum pp 27-28]{SibleyWentworth}).

Let $\mathcal{F}$ be a coherent sheaf over $X$, and $\mathcal{J} = \sqrt{\text{Ann}(\mathcal{F})}$ be the radical of its annihilator. We define $\text{gr}(\mathcal{F})$ as
\[\text{gr}(\mathcal{F}) = \bigoplus_{k \geq 0}\mathcal{J}^k\mathcal{F} / \mathcal{J}^{k+1}\mathcal{F}.\]

Note that this sum is clearly finite. 

\hfill

\begin{definition}\label{AlgMultiplicities}(see \cite[Erratum on pages 27-28]{SibleyWentworth})
	Let $\text{Gr}^S_{\o_0}(E)$ be an associated graded sheaf, and $Z_{alg}$ be the analytic set defined as in \ref{ZalgDefinition}. Let $\iota: Z \hookrightarrow Z_{alg}$ be an irreducible component of $Z_{alg}$ of codimension 2. The {\it algebraic multiplicity} $n_Z$ of $Z$ is the generic rank of $\iota^*\text{gr}\Big({(\text{Gr}^S_{\o_0}(E))^{**} \over \text{Gr}^S_{\o_0}(E)}\Big)$.	
\end{definition} 

\hfill

This definition appears in \cite[Errata, pp.27-28]{SibleyWentworth} as a replacement of \cite[Definition 2.4]{SibleyWentworth}.

\hfill

\begin{theorem}\label{SibleyWentworthTheorem}
	The subset $Z_{an}$ coincides with $Z_{alg}$ as sets. Moreover, for any codimension two irreducible subvariety $Z_\alpha \subset Z_{an}$, we have $Z_{\alpha} \subset Z_{alg}$, and its algebraic multiplicity coincides with $m_\alpha$.
\end{theorem}

\hfill

\begin{proof}
	The proof is essentially the same, as in \cite{SibleyWentworth}.  We just outline the main differences here.
	
	\textbf{Step 1:} The slicing lemma \cite[Lemma 4.1]{SibleyWentworth} ensures that for an embedded 4-disc $\Sigma \subset X$, which is a {\it transverse slice} for $Z_\alpha \subset Z_{b, an}$ (i.e. $\Sigma$ intersects transversely $Z_\alpha$ at one point), one has 
	\[m_\alpha = \lim_{j \to + \infty}\int_\Sigma c_2(A_j) - c_2(A_\infty) = \lim_{j \to + \infty}{1 \over 8\pi^2}\int_\Sigma \Tr (F_{A_j} \wedge F_{A_j}) - \Tr (F_{A_\infty} \wedge F_{A_\infty}) .\]
	The cancellation of $\Tr F_{A_j} \wedge \Tr F_{A_j} - \Tr F_{A_\infty} \wedge \Tr F_{A_\infty}$ is justified in the proof of \ref{ZanIsAnalytic}, and the proof of the equality above follows from the same argument, which Sibley and Wentworth use to show that the integral $\int_{\p \Sigma }\text{CS}(A_j, A_\infty)$ tends to zero as $j \to + \infty$. 
	
	\textbf{Step 2:}
	Next, \cite[Lemma 4.3]{SibleyWentworth} is a general complex-analytic fact, and does not require any information about our connections $A_j$. Also, we can use \cite[Theorem 1.3]{SibleyWentworth}. This theorem says, that for a holomorphic vector bundle $E$ with metric $H$, the following holds:
	\[\text{ch}_2(\text{Gr}(E),\tilde{H}) - \text{ch}_2(E, H) = \sum_\alpha n_\alpha W_\alpha + dd^c\Psi.\]
	
	Here $\tilde{H}$ is an admissible metric on the free part of $\text{Gr}^S(E)^{**}$ (i.e. with $L^2$-curvature and with bounded $|\Lambda_{\o_0}F_{\tilde{H}}|$; see \cite{BandoSiu}), $W_\alpha$ is a codimension 2 component of $Z_{alg}$, and $n_\alpha$ are defined in \ref{AlgMultiplicities}.
	
	\textbf{Step 3:} Finally, the proof of \cite[Proposition 4.2]{SibleyWentworth} carries over verbatim. Let $Z_\alpha \subset Z_{alg}$ with algebraic multiplicity $n_\alpha$. Upon taking a modification of $X$ (see \cite[Lemma 4.3]{SibleyWentworth}), we can take a closed 4-dimensional submanifold $\Sigma \subset X$, such that $[\Sigma]\cdot[Z_\alpha] \neq 0$, and $\Sigma$ intersects $Z_\alpha$ transversely at finitely many points $z_1, \dots, z_k$, and avoids other components. Then, using Step 2, we can repeat arguments of Sibley and Wentworth, and show that
	\begin{align*}
		&n_\alpha[\Sigma]\cdot[Z_\alpha] = {1 \over 8\pi^2}\int_{\Sigma \cap (\cup_{a=1}^k B_{\ve}(z_a))} \Tr (F_{A_j} \wedge F_{A_j}) - \Tr (F_{A_\infty} \wedge F_{A_\infty})
		\\
		&-\int_{\Sigma \cap \p (\cup_{a=1}^k B_{\ve}(z_a))} \text{CS}(A_j, A_\infty).
	\end{align*}
	
	By Step 1, we can take the limit of the right hand side, and obtain that it is equal to $m_\alpha[\Sigma]\cdot[Z_\alpha]$. 
	
	Consider the current 
	\[\sum_{\alpha}m_\alpha Z_\alpha - \sum_\beta n_\beta W_\beta.\]
	Here $Z_\alpha \subset Z_{an}$, and $W_\beta \subset Z_{alg}$ are codimension 2 irreducible components of $Z_{an}$ and $Z_{alg}$. As in \cite{SibleyWentworth} we see, that this current is nonnegative and exact. Hence it is zero, and any codimension 2 component of $Z_{an}$ lies in $Z_{alg}$ with the same multiplicity. 
	
	\textbf{Step 4:} We showed, that the codimension 2 components of $Z_{an}$ and $Z_{alg}$ coincide, and their multiplicities are the same. Thus, the remaining part is to show that the higher codimension components coincide as sets. This can be done again as Sibley and Wentworth did after the proof of \cite[Proposition 4.2]{SibleyWentworth}.
\end{proof}

\hfill

\hfill

\hfill

\hfill

\noindent Beijing Institute of Mathematical Sciences and Applications, China

\noindent nklemyatin@bimsa.cn 
%

\end{document}